\documentclass{amsart}

\usepackage{amssymb,latexsym}
\usepackage[latin1]{inputenc} 
\usepackage[dvips]{graphicx,color,psfrag}

\numberwithin{equation}{section}
\numberwithin{figure}{section}
\numberwithin{table}{section}

\theoremstyle{plain}
\newtheorem{prop}{Proposition}[section]
\newtheorem{lemma}[prop]{Lemma}
\newtheorem{theorem}[prop]{Theorem}

\newtheorem{conjecture}[prop]{Conjecture}

\theoremstyle{definition}
\newtheorem{defin}{Definition}[section]

\newenvironment{proofoft}[1]{{\em Proof of Theorem #1. }}{$\Box$ \vspace{1em}}

\newcommand{\Z}{{\mathbb Z}}
\newcommand{\Zd}{{\mathbb Z}^d}

\newcommand{\Remarks}{\em Remarks: \rm}

\def\P{{\mathbf P}}

\begin{document}

\title[Spin systems in a randomly evolving environment]
{Attractive nearest-neighbor spin systems on the integers in a randomly evolving environment}
\date{\today}

\author[Marcus Warfheimer]{Marcus Warfheimer}
\address{Department of Mathematical Sciences, Chalmers University of Technology and University of Gothenburg, SE-41296 Gothenburg, Sweden}
\email{marcus.warfheimer@gmail.com}
\urladdr{http://www.math.chalmers.se/~warfheim}
\thanks{Research partially supported by the Göran Gustafsson Foundation for Research in Natural Sciences and Medicine.}


\begin{abstract}
We consider spin systems on $\Z$ (i.e.\ interacting particle systems on $\Z$ in which each coordinate only has two possible values and only one coordinate changes in each transition) whose rates are determined by another process, called a background process. A canonical example is the contact process in randomly evolving environment, introduced and analysed by Broman and further studied by Steif and the author, where the marginals of the background process independently evolve as 2-state Markov chains and determine the recovery rates for a contact process. We prove that, if the background process has a unique stationary distribution and if the rates satisfy a certain positivity condition, then there are at most two extremal stationary distributions. The proof follows closely the ideas of Liggett's proof of a corresponding  theorem for spin systems on $\Z$ without a background process. 
\end{abstract}

\maketitle

\section{Introduction}

The contact process in a random environment, in which the rates are taken to be random variables and then fixed in time, has been studied the last twenty years, see for example \cite{Durrett1,Extinction1,Liggett1,Survival1}. However, recently Broman \cite{Broman} introduced a variant where the environment changes in time in a Markovian way. (See also \cite{JSMW} for further analysis concerning that process.) More precisely, he considered the Markov process $\{(B_t,C_t)\}_{t\geq 0}$ on $\{0,1\}^{\Zd}\times\{0,1\}^{\Zd}$ described by the following rates at a site $x$:
\[
\begin{array}{ll}
\textrm{transition} \qquad & \qquad \textrm{rate} \vspace{7pt} \\ 
(0,0) \rightarrow (0,1) \qquad &  \qquad\sum_{y\sim x} C(y) \\
(1,0) \rightarrow (1,1) \qquad &  \qquad\sum_{y\sim x} C(y) \\
(0,1) \rightarrow (0,0) \qquad &  \qquad\delta_0 \\
(1,1) \rightarrow (1,0) \qquad &  \qquad\delta_1 \\
(0,0) \rightarrow (1,0) \qquad &  \qquad\gamma p \\
(0,1) \rightarrow (1,1) \qquad &  \qquad\gamma p \\
(1,0) \rightarrow (0,0) \qquad &  \qquad\gamma (1-p) \\
(1,1) \rightarrow (0,1) \qquad &  \qquad\gamma (1-p) 
\end{array}
\]
where $\gamma, \delta_0, \delta_1>0$ with $ \delta_1 \leq\delta_0 $ and $p \in [0,1]$. In other words, the background process evolves independently for each site and determines the recovery rate for the right marginal in the following way: At a given site $x$ and time $t$, the rate is $\delta_0$ or $\delta_1$ depending on whether $B_t(x)=0$ or $B_t(x)=1$. Broman called $\{(B_t,C_t)\}$ the contact process in a randomly evolving environment, abbreviated CPREE. In this paper we study processes in one dimension with the same structure: a background process influencing another interacting particle system, but here both processes are more general. We prove, under certain conditions on the rates, that we have at most two extremal invariant distributions.

\section{The model and main result}\label{paper2:model}

We consider the Markov process, $\{(\beta_t,\eta_t)\}_{t \geq 0}$ on $\{0,1\}^{\Z}\times\{0,1\}^{\Z}$ described by the following rates at a site $x$:

\[
\begin{array}{lll}
\text{transition} \qquad  \qquad &\text{rate}  \vspace{7pt} \\ 
(\beta,\eta) \rightarrow (\beta,\eta_x) \qquad &  c_0(x,\eta) & \textrm{if $\:$ $\beta(x)=0$} \\
(\beta,\eta) \rightarrow (\beta,\eta_x) \qquad &  c_1(x,\eta) & \textrm{if $\:$ $\beta(x)=1$} \\
(\beta,\eta) \rightarrow (\beta_x,\eta) \qquad &  b(x,\beta) 
\end{array}
\]
Here $c_0(x,\eta)$, $c_1(x,\eta)$ and $b(x,\beta)$ are given rate functions where the first two satisfy
\begin{equation}\label{paper2:attr1}
\begin{split}
&c_0(x,\eta)\leq c_1(x,\eta) \quad \text{if} \quad \eta(x)=0, \\
&c_1(x,\eta)\leq c_0(x,\eta) \quad \text{if} \quad \eta(x)=1, 
\end{split}
\end{equation}
and all three satisfy the following attractivity condition:
\begin{defin}\label{paper2:attr}
A spin system on $\{0,1\}^{\Z}$, with rates $c(x,\eta)$ is said to be attractive if whenever $\eta \leq \eta^\prime$,
\begin{equation}\label{paper2:attr2}
\begin{split}
&c(x,\eta)\leq c(x,\eta^\prime) \quad \text{if} \quad \eta(x)=\eta^\prime(x)=0, \\
&c(x,\eta)\geq c(x,\eta^\prime) \quad \text{if} \quad \eta(x)=\eta^\prime(x)=1.  
\end{split}
\end{equation}
\end{defin}
Here, $\leq$ refers to the usual partial ordering on $\{0,1\}^{\Z}$, i.e., $\eta\leq\eta^\prime$ if and only if $\eta(x)\leq\eta^\prime(x)$ for all $x \in \Z$. We also assume that the rate functions are translation invariant and that the rates $c_0(x,\eta)$, $c_1(x,\eta)$ only depend on $\eta$ through 
\[
\{\eta(x-1),\eta(x),\eta(x+1)\}.
\]
Moreover, to ensure that we have a well defined process we will assume that 
\[
\displaystyle\sum_{y\in\Z}\sup_{\beta\in\{0,1\}^{\Z}}|b(0,\beta)-b(0,\beta_y)|<\infty.
\]
In other words, the rates for the system are completely described by $b(x,\beta)$ and the $16$ parameters determining $c_0$ and $c_1$. To describe the values we will use the following notation: 
\[
c_i(001)=c_i(x,\eta) \quad \text{when} \quad \eta(x-1)=0,\, \eta(x)=0 \;\, \text{and} \;\, \eta(x+1)=1. 
\]
We always refer to the left marginal as the \emph{background process}. Furthermore, note that we can equivalently view our process on $\{0,1\}^{\Z\times\{0,1\}}$ and that the conditions \eqref{paper2:attr1} and \eqref{paper2:attr2} then mean that the whole process is attractive on that space. (Definition~\ref{paper2:attr} can of course be generalized to $\{0,1\}^S$ where $S$ is countable.) The attractivity can be used to show (via monotonicity) the existence of two extremal stationary distributions $\nu_0$ and $\nu_1$ defined by
\begin{equation*}
\nu_0 = \lim_{t \to\infty}\delta_0 S(t) \quad \qquad \nu_1 = \lim_{t \to\infty}\delta_1 S(t), 
\end{equation*}
where $\delta_0$ and $\delta_1$ denote the point masses corresponding to the elements $\eta\equiv 0$ and $\eta\equiv 1$ in $\{0,1\}^{\Z\times\{0,1\}}$ and $\{S(t)\}_{t\geq 0}$ denotes the semigroup associated to $\{(\beta_t,\eta_t)\}_{t\geq 0}$. The main result here is that, if the background process has a unique stationary distribution and the rates $c_0$, $c_1$ satisfy a certain positivity condition, then $\nu_0$ and $\nu_1$ are the only extremal stationary distributions. Let $\mathcal I$ denote the set of stationary distributions for the process and let $\mathcal I_e$ denote its extreme points. Furthermore, define 
\begin{align*}
C_1&=\{\, c_i(100)+c_j(110),\,c_i(001)+c_j(011),\\
&\quad c_i(011)+c_j(110),\,c_i(100)+c_j(001),\,i=0,1,\,j=0,1 \,\}
\end{align*}
and let
\[
C=\min\left(C_1 \right).
\]

Before we state our main result, we want to emphasize that the case with no background process has been studied before by Liggett. The proof of our main result follows closely the ideas of his proof. To state his result, let $c(x,\eta)$ be a rate function for an attractive, translation invariant, nearest-neighbor spin system $\{\eta_t\}_{t\geq 0}$ on $\{0,1\}^{\Z}$ and define $\mu_i=\lim_{t\to\infty}\delta_i T(t)$, $i=0,1$, where $\delta_i$ is the point mass corresponding to the element $\eta\equiv i$ in $\{0,1\}^{\Z}$ and $\{T(t)\}_{t\geq 0}$ denotes the semigroup associated to $\{\eta_t\}_{t\geq 0}$. Moreover, let $\mathcal J_e$ denote the extreme points of the set of stationary distributions for $\{\eta_t\}_{t\geq 0}$.
\begin{theorem}[Liggett]\label{paper2:Liggett}
Suppose  
\begin{equation}\label{paper2:Liggett_assump}
c(x,\eta)+c(x,\eta_x)>0 \quad \text{whenever} \quad \eta(x-1)\neq\eta(x+1).
\end{equation}
Then ${\mathcal J}_e=\{\mu_0, \mu_1\}$.  
\end{theorem}
For a proof, see \cite{Liggett_spin_systems} or \cite[p.~145-152]{Liggett85}. In fact, he also proved that if condition \eqref{paper2:Liggett_assump} fails, then $\mathcal J_e$ contains infinitely many points, see \cite[p.~145]{Liggett85}.
\begin{theorem}\label{paper2:mainthm}
Suppose that the background process has a unique stationary distribution and assume $C>0$. Then ${\mathcal I}_e=\{\nu_0, \nu_1\}$.  
\end{theorem}


\noindent
\Remarks 
\begin{itemize}
\item[(i)] From \cite[p.~152]{Liggett85} we get that Theorem~\ref{paper2:Liggett} is equivalent to the statement that \eqref{paper2:Liggett_assump} and 
\begin{align*}
c(011)+c(110)&>0 \\ 
c(100)+c(001)&>0
\end{align*}
implies ${\mathcal J}_e=\{\mu_0, \mu_1\}$. By letting $c=c_0=c_1$, it is now clear that Theorem~\ref{paper2:mainthm} covers Theorem~\ref{paper2:Liggett}. 
\item[(ii)] The hypotheses in Theorem~\ref{paper2:mainthm} are true for the CPREE described in the introduction. Indeed, if $c_0$ and $c_1$ satisfy \eqref{paper2:attr1} and are symmetric under reflections, i.e.\ 
\begin{align*}
c_i(100) & = c_i(001) \\
c_i(110) & = c_i(011), \quad i=0,1
\end{align*}
then $C>0$ if and only if $c_0(001)>0$ and $c_1(011)>0$.  
\item[(iii)] Note that we are not assuming independence or even nearest-neighbor interaction between coordinates in the background process.
\item[(iv)] To see that the conclusion may fail if we drop the assumption about a unique stationary distribution for the background process, let $b(x,\beta)$, in addition to being attractive and translation invariant, be nearest-neighbor with $b(000)=b(111)=0$ and satisfiy
\[
b(x,\beta)+b(x,\beta_x)>0 \quad \text{whenever}\quad \beta(x-1)\neq\beta(x+1).
\]
Let $c_0=c_1$ be the rates corresponding to a supercritical contact process on $\Z$. Then
\[
{\mathcal I}_e=\{\,\delta_0\times\delta_0,\,\delta_0\times\bar{\nu},\,\delta_1\times\delta_0,\,\delta_1\times\bar{\nu} \,\},
\]
where $\delta_0$, $\delta_1$ are the point masses corresponding to the elements $\eta\equiv 0$ and $\eta\equiv 1$ in $\{0,1\}^{\Z}$ respectively and $\bar{\nu}$ denotes the upper invariant measure for the contact process. 
\item[(v)] If we take the same background process, but instead let $c_0=c_1$ be the rates for a subcritical contact process, we see that the condition about a unique stationary distribution for the background process  is not necessary for having only two extremal stationary distributions.
\item[(vi)] To see that the conclusion may fail if $C=0$, let $b(x,\beta)$ be a rate function such that $\{\beta_t\}_{t\geq 0}$ has the point mass at $\beta\equiv 1$ as its unique stationary distribution and let $c_1$ satisify
\[
c_1(001)+c_1(011)=0.
\]
It is easy to check that for each $n\in \Z$, $\delta_1\times\delta_{\eta^n}$ is an extremal stationary distribution where
\[
\eta^n(x)=\begin{cases}
1 &\text{if $x\geq n$} \\
0 &\text{if $x<n$}.
\end{cases} 
\]   
\end{itemize}
A natural next step is to ask when there is a unique stationary distribution, i.e.\ when $\nu_0=\nu_1$. In the case of no background process, Gray proved in \cite{Gray_pos_rates} that there can only be one stationary distribution provided that the rates are strictly positive. We conjecture an analogous statement in our situation.  
\begin{theorem}[Gray]
If $c(x,\eta)>0$ for all $x\in\Z$ and $\eta\in\{0,1\}^\Z$, then $\mu_0=\mu_1$.
\end{theorem}
\begin{conjecture}\label{paper2:conj_pos_rates}
Suppose that the background process has a unique stationary distribution and assume that $c_i(x,\eta)>0$ for all $x$, $\eta$, $i=1,2$. Then $\nu_0=\nu_1$. 
\end{conjecture}

The rest of the paper is organized as follows. In Section~\ref{paper2:proof_mainthm} we prove Theorem~\ref{paper2:mainthm} and in Section~\ref{paper2:conj} we discuss Conjecture~\ref{paper2:conj_pos_rates}.

\section{Proof of Theorem~\ref{paper2:mainthm}}\label{paper2:proof_mainthm}

In the proof, we  make extensive use of a maximal type coupling which we now describe. Denote
\begin{align*}
U &=\{0,1\}^{\Z}, & V=\{\,(\eta,\gamma,\xi)\in U^3:\,\eta\leq\gamma\leq\xi\, \} \quad \textrm{and}\quad & W=U\times V.
\end{align*}
The coupled process $(\beta_t, \eta_t,\gamma_t,\xi_t)$, which we now define, lives on $W$ and its flip rates are described as follows: First, let flips of the type
$$
(\beta,\eta,\gamma,\xi) \rightarrow (\beta_x,\eta,\gamma,\xi)
$$
occur at rate $b(x,\beta)$. 
\begin{table}[!h]
\begin{center}
\begin{tabular}{c|c|c|c|c|}
& (0,0,0,0) & (0,0,0,1) & (0,0,1,1) & (0,1,1,1) \\
\hline
(0,0,0,0) & -- & $c_0(x,\xi)-c_0(x,\gamma)$ & $c_0(x,\gamma)-c_0(x,\eta)$ & $c_0(x,\eta)$ \\
\hline
(0,0,0,1) & $c_0(x,\xi)$ & -- & $c_0(x,\gamma)-c_0(x,\eta)$ & $c_0(x,\eta)$ \\
\hline
(0,0,1,1) & $c_0(x,\xi)$ & $c_0(x,\gamma)-c_0(x,\xi)$ & -- & $c_0(x,\eta)$ \\
\hline
(0,1,1,1) & $c_0(x,\xi)$ & $c_0(x,\gamma)-c_0(x,\xi)$ & $c_0(x,\eta)-c_0(x,\gamma)$ & -- \\
\hline 
\end{tabular} 
\caption{Transition rates when the background process is in state 0.}
\label{paper2:rates1}
\end{center}
\end{table}
\begin{table}[!h]
\begin{center}
\begin{tabular}{c|c|c|c|c|}
& (1,0,0,0) & (1,0,0,1) & (1,0,1,1) & (1,1,1,1) \\
\hline
(1,0,0,0) & -- & $c_1(x,\xi)-c_1(x,\gamma)$ & $c_1(x,\gamma)-c_1(x,\eta)$ & $c_1(x,\eta)$ \\
\hline
(1,0,0,1) & $c_1(x,\xi)$ & -- & $c_1(x,\gamma)-c_1(x,\eta)$ & $c_1(x,\eta)$ \\
\hline
(1,0,1,1) & $c_1(x,\xi)$ & $c_1(x,\gamma)-c_1(x,\xi)$ & -- & $c_1(x,\eta)$ \\
\hline
(1,1,1,1) & $c_1(x,\xi)$ & $c_1(x,\gamma)-c_1(x,\xi)$ & $c_1(x,\eta)-c_1(x,\gamma)$ & -- \\
\hline 
\end{tabular} 
\caption{Transition rates when the background process is in state 1.}
\label{paper2:rates2}
\end{center}
\end{table}

Then, let the other three marginals flip according to Tables~\ref{paper2:rates1} and \ref{paper2:rates2}. These tables should be interpreted as follows. For example, when $\beta_t(x)=0$, $\eta_t(x)=0$, $\gamma_t(x)=0$ and $\xi_t(x)=1$, $\xi_t(x)$ will flip alone at rate $c_0(x,\xi_t)$, $\gamma_t(x)$ will flip alone at rate $c_0(x,\gamma_t)-c_0(x,\eta_t)$ and $\eta_t(x)$ and $\gamma_t(x)$ flip together at rate $c_0(x,\eta_t)$. Note that the pairs $\{(\beta_t,\eta_t)\}$, $\{(\beta_t,\gamma_t)\}$, $\{(\beta_t,\xi_t)\}$ each evolve as the original Markov process and that the second, third and fourth marginals try to flip together as much as possible. Also, observe that the background process is not allowed to flip together with any of the other processes. 

As in the proof of Theorem~\ref{paper2:Liggett}, the proof of Theorem~\ref{paper2:mainthm} consists of several lemma concerning certain functionals of the process. For $m\leq n$, let $f_{m,n}(\beta,\eta,\gamma,\xi)$ be the number of intervals of zeros and ones in $\gamma$ between $m$ and $n$ (including $m$ and $n$), counted only where $\eta$ and $\xi$ differ. Furthermore, let 
\[
m\leq x_1 < x_2<\ldots<x_k\leq n, 
\]
be all those $x$'s between $m$ and $n$ for which $\eta(x)=0$ and $\xi(x)=1$. For $l\geq 1$, define 
\begin{align*}
g_{m,n}^l(\beta,\eta,\gamma,\xi)&=\text{number of $i$ such that $i\geq 1$, $i+l+1\leq k$ and} \\
\gamma(x_i)&\neq \gamma(x_{i+1})=\gamma(x_{i+2})=\ldots=\gamma(x_{i+l})\neq\gamma(x_{i+l+1}).
\end{align*}
In other words, $g_{m,n}^l(\beta,\eta,\gamma,\xi)$ is the number of interior intervals of zeros and ones of length $l$ in $\gamma$ between $m$ and $n$, counted only where $\eta$ and $\xi$ differ. For example if, 
\[
\begin{array}{ccccccccccccc|c}
\cdots & 1 & 0 & 1 & 1 & 1 & 1 & 1 & 0 & 1 & 1 & 1 & \cdots & \xi \\
\cdots & 1 & 0 & 1 & 1 & 0 & 0 & 1 & 0 & 1 & 1 & 0 & \cdots & \gamma \\
\cdots & 1 & 0 & 0 & 0 & 0 & 0 & 0 & 0 & 0 & 0 & 0 & \cdots & \eta \\
\cdots & 1 & 0 & 1 & 0 & 0 & 1 & 1 & 1 & 0 & 1 & 1 & \cdots & \beta  \\ 
         & m &   &   &   &   &   &   &   &   &   & n              
\end{array}
\]
then $f_{m,n}=4$, $g_{m,n}^2=1$, $g_{m,n}^3=1$ and $g_{m,n}^l=0$ when $l\notin \{\,2,3\,\}$. Let 
\[
K=\max\left(\max_{\eta}c_0(x,\eta),\,\max_{\eta}c_1(x,\eta)\right)
\]
and denote the set of stationary distributions and the generator of the coupled process by $\tilde{\mathcal I}$ and $\tilde{\Omega}$ respectively. Furthermore, for a given set $\mathcal A$, denote the set of extreme points by $\mathcal A_e$. The first lemma concerns certain basic properties of $f_{m,n}$ and $g_{m,n}^l$.
\begin{lemma}\label{paper2:basic_prop}
\begin{align*}
\text{a) } & f_{m,n}, g_{m,n}^l \text{ are increasing when $n$ increases or $m$ decreases.} \\
\text{b) } & f_{m,n}\leq 2+ \sum_{l=1}^{\infty}g_{m,n}^l. \\
\text{c) } & \sum_{l=1}^{\infty} l g_{m,n}^l \leq n-m+1. \\
\text{If } & \nu \in \tilde{\mathcal I}, \\
\text{d) } & C \int g_{m,n}^1 \,d\nu \leq K \int \left[f_{m-1,n}+f_{m,n+1}-2 f_{m,n}\right]\, d\nu, \text{ for $m\leq n$}\\
\text{e) } & C \int g_{m,n}^{l+1}\, d\nu \leq 12\text{Kl} \int g_{m,n}^{l}\, d\nu,\text{ for $m\leq n$, $l\geq 1$}.
\end{align*}
\end{lemma}
\begin{proof}
a), b) and c) follow directly from the definitions. For d) and e) assume $\nu \in \tilde{\mathcal I}$. Note that $f_{m,n}$ and $g_{m,n}^l$ are cylinder functions so that
\begin{align}\label{paper2:cyl}
\int \tilde{\Omega} f_{m,n} \,d \nu = \int \tilde{\Omega} g_{m,n}^l \,d \nu=0.
\end{align}
For cylinder function $f$, the generator has the form 
\begin{equation}\label{paper2:gen}
\begin{split}
\tilde{\Omega} f(\beta,\eta,\gamma,\xi)&=\sum_{(\beta,\bar{\eta},\bar{\gamma},\bar{\xi})} c(\beta,\eta,\gamma,\xi,\bar{\eta},\bar{\gamma},\bar{\xi})\left(f(\beta,\bar{\eta},\bar{\gamma},\bar{\xi})- f(\beta,\eta,\gamma,\xi) \right) \\
&\quad+\sum_{x} b(x,\beta) \left( f(\beta_x,\eta,\gamma,\xi)- f(\beta,\eta,\gamma,\xi) \right)
\end{split}
\end{equation}
where the first sum is over all possible transitions when the second, third or fourth marginal flip. (Recall that the first marginal is not allowed to flip together with any of the others.) Here, since both $f_{m,n}$ and $g_{m,n}^l$ do not depend on $\beta$, the second sum is zero, so our task is to calculate the first part. For this, we follow the approach in \cite[Lemma 3.7]{Liggett85}. The argument given here is almost the same as in \cite{Liggett85}, we supply it for the sake of completeness. Let $(\beta,\eta,\gamma,\xi)$ be fixed and note that the only way $f_{m,n}$ can increase because of a flip is if $f_{m-1,n}=f_{m,n}+1$ or $f_{m,n+1}=f_{m,n}+1$. In the first case the flip must occur at $x=m$ and in the second at $x=n$. The rate for such a flip is at most $K$ so the positive terms in \eqref{paper2:gen} are bounded above by 
\[
K\left[f_{m-1,n}+f_{m,n+1}-2f_{m,n}\right].
\]
Furthermore, there are $g_{m,n}^1$ sites $x$ where a flip decreases $f_{m,n}$ by two. At such an $x$, $\gamma(x)=0$ or $\gamma(x)=1$. Assume $\gamma(x)=1$. Then we necessarely have $\gamma(x-1)=\eta(x-1)$ and $\gamma(x+1)=\eta(x+1)$. Therefore, the flip rate at $x$ becomes
\[
c_0(x,\gamma)+c_0(x,\eta)=
\begin{cases}
c_0(010)+c_0(000) & \text{if $\gamma(x-1)=0$, $\gamma(x+1)=0$}, \\
c_0(011)+c_0(001) & \text{if $\gamma(x-1)=0$, $\gamma(x+1)=1$}, \\
c_0(110)+c_0(100) & \text{if $\gamma(x-1)=1$, $\gamma(x+1)=0$}, \\
c_0(111)+c_0(101) & \text{if $\gamma(x-1)=1$, $\gamma(x+1)=1$},
\end{cases} 
\]
when $\beta(x)=0$ and
\[
c_1(x,\gamma)+c_1(x,\eta)=
\begin{cases}
c_1(010)+c_1(000) & \text{if $\gamma(x-1)=0$, $\gamma(x+1)=0$}, \\
c_1(011)+c_1(001) & \text{if $\gamma(x-1)=0$, $\gamma(x+1)=1$}, \\
c_1(110)+c_1(100) & \text{if $\gamma(x-1)=1$, $\gamma(x+1)=0$}, \\
c_1(111)+c_1(101) & \text{if $\gamma(x-1)=1$, $\gamma(x+1)=1$},
\end{cases} 
\]
when $\beta(x)=1$. Also the attractivity condition gives
\begin{align*}
c_i(010) & \geq \max\{\,c_i(011),\,c_i(110)\,\} \\
c_i(101) & \geq \max\{\,c_i(001),\,c_i(100)\,\}, \quad i=0,1
\end{align*}
and so the rates above are bounded below by $C/2$. The same argument works if $\gamma(x)=0$ and so we can conclude that the negative terms in \eqref{paper2:gen} are bounded above by $-C g_{m,n}^1$. We get the estimate
\[
\tilde{\Omega} f_{m,n} \leq K\left[f_{m-1,n}+f_{m,n+1}-2f_{m,n} \right] - C g_{m,n}^1
\]
which via \eqref{paper2:cyl} gives d). For e), note that $g_{m,n}^l$ can only decrease via flips at no more than $l g_{m,n}^l$ sites or their neighbors, i.e.\ in total at most $3l g_{m,n}^l$ sites. The rate for such a flip is bounded by $2K$ and $g_{m,n}^l$ can at most decrease by two. The negative terms in the generator are therefore bounded below by $-12Kl g_{m,n}^l$. Furthermore, $g_{m,n}^l$ can increase at no fewer than $g_{m,n}^{l+1}$ pair of sites. These pair of sites are the endpoints of an interval of length $l+1$. To get a lower bound on the flip rate for such endpoints, let $x<y$ denote such a pair and suppose $\gamma(x)=\gamma(y)=1$. Then we have $\gamma(x-1)=\eta(x-1)$ and $\gamma(y+1)=\eta(y+1)$. The flip rate at $x$ is at least $c_i(100)$ if $\gamma(x-1)=\eta(x-1)=1$, $\beta(x)=i$ and at least $c_i(011)$ if $\gamma(x-1)=\eta(x-1)=0$, $\beta(x)=i$. In a similar fashion, the flip rate at $y$ is at least $c_i(001)$ if $\gamma(y+1)=\eta(y+1)=1$, $\beta(y)=i$ and at least $c_i(110)$ if $\gamma(y+1)=\eta(y+1)=0$, $\beta(y)=i$. In either case the sum of the flip rates for the pair is always at least $C$. The same statement holds if $\gamma(x)=\gamma(y)=0$ and so we obtain that the positive terms in the generator expression are bounded below by $C g_{m,n}^{l+1}$. Hence, we get the estimate
\[
\tilde{\Omega} g_{m,n} \geq C g_{m,n}^{l+1} - 12Kl g_{m,n}^l.
\]
Equation \eqref{paper2:cyl} then finally gives us
\[
C \int g_{m,n}^{l+1} \,d\nu \leq 12Kl \int g_{m,n}^{l} \,d\nu
\]
and the proof is complete.     
\end{proof} 

\noindent
Denote
\begin{align*}
A_1 & = \{\,(\beta,\eta,\gamma,\xi)\in W:\,\gamma\equiv\eta \,\}, \\
A_2 & = \{\,(\beta,\eta,\gamma,\xi)\in W:\,\gamma\equiv\xi \,\}, \\
\begin{split}
A_3 & = \{\,(\beta,\eta,\gamma,\xi)\in W\setminus A_1\cup A_2:\,\exists x \in \Z \text{ such that } \\
&\qquad \gamma(y)=\eta(y) \text{ when } y \leq x \text{ and } \gamma(y)=\xi(y) \text{ when } y>x\,\}, \end{split} \\
\begin{split}
A_4 & = \{\,(\beta,\eta,\gamma,\xi)\in W\setminus A_1\cup A_2:\,\exists x \in \Z \text{ such that } \\
&\qquad\gamma(y)=\xi(y) \text{ when } y \leq x \text{ and } \gamma(y)=\eta(y) \text{ when } y>x\,\}, \end{split}
\end{align*}
\begin{lemma}\label{paper2:conc}
Assume $C>0$. Then
\begin{align*}
\text{a) }& \nu\in \tilde{\mathcal I} \quad \implies \quad \nu\left(A_1\cup A_2\cup A_3\cup A_4 \right)=1, \\
\text{b) }& \nu\in \tilde{\mathcal I}_e \quad \implies \quad \nu\left(A_i\right)=1 \,\,\text{ for some $i$.} 
\end{align*}
\end{lemma}
\begin{proof}
$b)$ follows from a) since $A_i$ is closed for the coupled process in the sense that
\[
\P^{(\beta,\eta,\gamma,\xi)}[\,(\beta_t,\eta_t,\gamma_t,\xi_t)\in A_i \,]=1 \quad \forall t>0
\]
whenever $(\beta,\eta,\gamma,\xi)\in A_i$. To prove $a)$, suppose $\nu \in \tilde{\mathcal I}$. 
Since 
\[
\displaystyle\bigcup_{i=1}^4 A_i = \{\, g_{m,n}^l=0 \: \forall m \leq n,\, l \geq 1 \,\}
\]
we obtain that
\begin{equation}\label{paper2:int1}
\int g_{m,n}^l d\nu=0\text{ for all $m\leq n$, $l\geq 1$}
\end{equation}
is equivalent to 
\[
\nu\left(A_1\cup A_2\cup A_3\cup A_4 \right)=1.
\]
To see that \eqref{paper2:int1} holds, we proceed as in \cite[Lemma 3.10]{Liggett85}. Note that 
\begin{equation*}
f_{m-1,n} \leq f_{m,n}+1 \qquad  \text{and} \qquad  f_{m,n+1} \leq f_{m,n}+1 
\end{equation*}
and so parts $d)$ and $e)$ of Lemma~\ref{paper2:basic_prop} gives us
\begin{equation}\label{paper2:M}
M=\sup_{m\leq n} \int g_{m,n}^l \,d \nu < \infty, \quad \forall l \geq 1.
\end{equation}
Let $L \geq 1$. From part $b)$ of the same lemma, we get
\[
\frac{1}{n-m} \int f_{m,n} \,d \nu \leq \frac{2}{n-m} + \frac{1}{n-m} \int \sum_{l\geq 1} g_{m,n}^l \,d \nu.
\] 
Split the sum and now use part $c)$ of the lemma together with \eqref{paper2:M} to obtain that for any $L$
\[
\frac{1}{n-m} \int f_{m,n} \,d \nu \leq \frac{2}{n-m}+\frac{ML}{n-m}+\frac{1}{L}\left(1+ \frac{1}{n-m}\right),
\]
and so
\[
\limsup_{n-m \to \infty} \frac{1}{n-m} \int f_{m,n} \,d \nu \leq\frac{1}{L}.
\]
Since $L \geq 1$ was arbitrary we can conclude
\begin{equation}\label{paper2:lim0}
\lim_{n-m \to \infty} \frac{1}{n-m} \int f_{m,n} \,d \nu=0.
\end{equation}
Now, for $N \geq 1$, part $d)$ of Lemma~\ref{paper2:basic_prop} gives us
\begin{equation}\label{paper2:sum}
\begin{split}
&C\sum_{m=-N+1}^0 \sum_{n=0}^{N-1} \int g_{m,n}^1 \,d\nu \\
&\quad\leq K \sum_{m=-N+1}^0 \sum_{n=0}^{N-1} \int \left[f_{m-1,n}+f_{m,n+1}-2f_{m,n} \right]\,d\nu.
\end{split}
\end{equation}
After some cancellations in the sum to the right, we get
\[
\begin{split}
&\sum_{m=-N+1}^0 \sum_{n=0}^{N-1} \int \left[f_{m-1,n}+f_{m,n+1}-2f_{m,n} \right]\,d\nu \\ 
&\qquad\leq   \sum_{m=-N+1}^0 \int f_{m,N}\,d\nu+\sum_{n=0}^{N-1} \int f_{-N,n}\,d\nu
\end{split}
\]
and together with \eqref{paper2:lim0} and \eqref{paper2:sum} we obtain
\[
\lim_{N\to\infty}\frac{1}{N^2} \sum_{m=-N+1}^0 \sum_{n=0}^{N-1} \int g_{m,n}^1 \,d\nu=0.
\]
Using the monotonicity property of $g_{m,n}^1$ this implies $\int g_{m,n}^1 \,d\nu=0$ for all $m\leq n$ and part $e)$ of the lemma gives $\int g_{m,n}^l \,d\nu=0$ for all $l \geq 1$ and we are done with the proof.
\end{proof} 
We are soon ready for the proof of Theorem~\ref{paper2:mainthm}. However, in the proof we make use of a $5$-variant coupling $\{(\beta_t,\eta_t,\gamma_{1,t},\gamma_{2,t},\xi_t)\}$ of the one used so far. This coupling is also of maximal type and evolves on 
\[
X=\left\{\,(\beta,\eta,\gamma_1,\gamma_2,\xi)\in U^5:\:\eta \leq \gamma_1 \leq \xi, \,\eta \leq \gamma_2 \leq \xi\right\}
\]
in a way such that $\{(\beta_t,\eta_t,\gamma_{1,t},\xi_t)\}$ and $\{(\beta_t,\eta_t,\gamma_{2,t},\xi_t)\}$ evolve exactly as the previous described coupling. We can therefore apply all we have done so far to each of these processes. The last tool we need is to have existence of an extremal stationary distribution for the 5-variant coupled process, given extremal stationary distributions for the $\{(\beta_t,\eta_t)\}$ process. For a stochastic variable $X$ and a distribution $\mu$, let $X\sim\mu$ denote that $X$ is distributed according to $\mu$. Also, let $\mathcal I^5$ denote the set of stationary distributions for the $5$-variant coupled process on $X$. 
\begin{lemma}\label{paper2:statcoup}
Given  $\mu$, $\mu^\prime\in \mathcal I_e$ there exists $\nu((\beta,\eta,\gamma_1,\gamma_2,\xi)\in \cdot ) \in {\mathcal I^5_e}$ such that $(\beta,\eta) \sim \nu_0$, $(\beta,\gamma_1) \sim \mu$, $(\beta,\gamma_2) \sim \mu^\prime$ and $(\beta,\xi)\sim\nu_1$.
\end{lemma}
\begin{proof}
For any measure $\mu$ let $\mu_{ij}$ denote the projection to the \it i\rm th and \it j\rm th coordinate. Construct a coupling on $(\{0,1\}^{\Z}\times\{0,1\}^{\Z})^4$ of four $\{\beta_t,\eta_t\}$-processes such that the background processes agree as much as possible as well as the right marginals. Note that our $5$-variant coupling above can be identified with such a coupling started with all the background processes equal. 
Starting the coupling with
\[
\delta_{(\emptyset,\emptyset)}\times\mu\times\mu^\prime\times\delta_{(\Z,\Z)}
\]
and taking a suitable subsequence of Cesaro averages gives us a stationary distribution $\rho$ for the coupling and by projecting to the first, second, fourth, sixth and eighth coordinate we get a probability measure $\tilde\nu \in\mathcal I^5$ with
\[
\tilde{\nu}((\beta,\eta,\gamma_1,\gamma_2,\xi)\in U^5:\:\eta \leq \gamma_1 \leq \xi, \,\eta \leq \gamma_2 \leq \xi)=1.
\]
Here it is important to note that the set
\[
\begin{split}
&\{\,(\beta_1,\eta,\beta_2,\gamma_1,\beta_3,\gamma_2,\beta_4,\xi)\in U^8:\:\beta_1 \leq \beta_2 \leq \beta_4, \, \beta_1 \leq \beta_3 \leq \beta_4, \\
&\qquad\eta \leq \gamma_1 \leq \xi, \,\eta \leq \gamma_2 \leq \xi\}
\end{split}
\]
is closed under the evolution of the coupling and that the first, third, fifth and seventh coordinate are equal under $\rho$. Furthermore, it is clear that $\tilde\nu$ satisfies
\[
\tilde\nu_{12}=\nu_0, \quad \tilde\nu_{13}=\mu \quad \tilde\nu_{14}=\mu^\prime \quad \text{and} \quad \tilde\nu_{15}=\nu_1.
\]
Define 
\[
\mathcal B=\{\,\nu\in \mathcal I^5:\,\nu_{12}=\nu_0,\,\nu_{13}=\mu,\,\nu_{14}=\mu^\prime, \,\nu_{15}=\nu_1\,\}.
\]
$\mathcal B$ is non-empty by the above and is compact and convex. Hence, by the Krein-Milman theorem, $\mathcal B$ can be written as the closed convex hull of its extreme points. Therefore, since $\mathcal B\neq\emptyset$, we have ${\mathcal B}_e\neq\emptyset$. Hence, the proof is complete if ${\mathcal B}_e\subset{\mathcal I}^{5}_e$. Assume $\nu\in{\mathcal B}_e$ and let $\nu = \alpha\rho +(1-\alpha)\sigma$, where $0<\alpha<1$ and $\rho,\sigma\in\mathcal I^5$. If $\rho$, $\sigma\in \mathcal B$ we get $\nu=\rho=\sigma$ and we are done. In order to see this, let $(i,j)$ be one of the pairs $(1,2)$, $(1,3)$, $(1,4)$ or $(1,5)$. Since $\nu_{ij}=\alpha\rho_{ij} +(1-\alpha)\sigma_{ij}$, where $\rho_{ij},\sigma_{ij}\in \mathcal I$, and the left hand side is an element of $\{\nu_0,\mu,\mu^\prime,\nu_1\}\subseteq\mathcal I_e$, we obtain 
\begin{align*}
\nu_0=\rho_{12}  &= \sigma_{12}&   \mu=\rho_{13} &= \sigma_{13}  \\
\mu^\prime=\rho_{14} &= \sigma_{14}&   \nu_1=\rho_{15} &= \sigma_{15}
\end{align*}
and so $\rho$, $\sigma\in \mathcal B$.      
\end{proof}
\noindent
\begin{proofoft}{\ref{paper2:mainthm}}
We follow the steps in \cite[Theroem 3.13]{Liggett85}. Let $\mu_1\in {\mathcal I}_e$. Since $\nu_0 \leq \mu \leq \nu_1$ for every stationary distribution $\mu$, we can assume $\nu_0\neq \nu_1$. Let $\mu_2=\mu_1 \circ \theta_x^{-1}$, where $\theta_x$ is a translation by $x\in\Z$. Since the dynamics are translation invariant and $\mu_1 \in \mathcal I_e$, we get that $\mu_2 \in \mathcal I_e$. Let $\rho$ be an extremal stationary distribution for the $5$-variant coupling mentioned above with
\begin{align*}
(\beta,\eta)  &\sim \nu_0&   (\beta,\gamma_1) &\sim \mu_1 \\
(\beta,\gamma_2)  &\sim \mu_2&   (\beta,\xi) &\sim \nu_1
\end{align*}
Such a measure exists by Lemma~\ref{paper2:statcoup}. Let $\rho_1$ and $\rho_2$ be the distributions obtained from the projections
\begin{align*}
(\beta,\eta,\gamma_1,\gamma_2,\xi) &\to (\beta,\eta,\gamma_1,\xi) \\
(\beta,\eta,\gamma_1,\gamma_2,\xi) &\to (\beta,\eta,\gamma_2,\xi)
\end{align*}
respectively. Since $\rho_1,\rho_2\in \tilde{\mathcal I}_e$, Lemma~\ref{paper2:conc} gives
\begin{equation*}
\rho_1(A_i)=1 \quad \text{some $1\leq i\leq 4$} \quad \text{and} \quad \rho_2(A_i)=1 \quad \text{some $1\leq i\leq 4$}.
\end{equation*}
However, $\gamma_1$ and $\gamma_2$ are just translations of each other so there is an $i$ such that $\rho_1(A_i)=\rho_2(A_i)=1$. It follows that 
\[
\rho\Big((\beta,\eta,\gamma_1,\gamma_2,\xi):\: \sum_x \vert \gamma_1(x)-\gamma_2(x)\vert < \infty \Big)=1.
\]
Also, $(\gamma_{1,t},\gamma_{2,t})$ has the property that
\[
\P^{(\gamma,\gamma)}[\gamma_{1,t}=\gamma_{2,t}]=1 \quad \text{and} \quad \P^{(\gamma_1,\gamma_2)}[\gamma_{1,t}=\gamma_{2,t}]>0
\] 
whenever $\displaystyle\sum_x \vert \gamma_1(x)-\gamma_2(x)\vert < \infty$ and so since $\rho$ is stationary, we must in fact have
\[
\rho\Big((\beta,\eta,\gamma_1,\gamma_2,\xi):\: \gamma_1=\gamma_2 \Big)=1.
\]
This implies $\mu_1=\mu_2$, i.e.\ $\mu_1$ is translation invariant. Therefore $i$ equals $1$ or $2$ (recall $\nu_0\neq \nu_1$). If $i=1$, $\mu_1(U\times (\cdot))=\nu_0(U\times (\cdot))$ and since the background process has a unique stationary distribution we must also have $\mu_1((\cdot)\times U)=\nu_0((\cdot)\times U)$. But since $\nu_0 \leq \mu_1$ this yields $\mu_1=\nu_0$. If $i=2$ we get in a similar way that $\mu_1=\nu_1$.  
\end{proofoft}

\section{Discussion of Conjecture~\ref{paper2:conj_pos_rates}}\label{paper2:conj}

We begin by describing a graphical representation which may be useful for a possible proof of Conjecture~\ref{paper2:conj_pos_rates}. The representation is similar as in \cite{Gray_pos_rates} and we will explain it in a quite informal way. For simplicity, we will assume that the rates for the background process, in addition to attractive and translation invariant, also are uniformly bounded. (Of course, our assumptions on $c_0$ and $c_1$ from Section~\ref{paper2:model} imply that they are also uniformly bounded.) For $x\in\Z$, define
\[
\begin{split}
\bar{b}_x&=\sup_{\beta:\,\beta(x)=0}b(x,\beta)+\sup_{\beta:\,\beta(x)=1}b(x,\beta) \\
\bar{c}_x^0&=\sup_{\eta:\,\eta(x)=0}c_0(x,\eta)+\sup_{\eta:\,\eta(x)=1}c_0(x,\eta) \\
\bar{c}_x^1&=\sup_{\eta:\,\eta(x)=0}c_1(x,\eta)+\sup_{\eta:\,\eta(x)=1}c_1(x,\eta) \\
\bar{c}_x&=\bar{c}_x^0+\bar{c}_x^1.
\end{split}
\]
Define the following collection of independent random variables on some probability space $(\Omega, \mathcal F, \P)$:
\begin{itemize}
\item[--] $B_j(x)$ exponentially distributed with mean $1/\bar{b}_x$, $j\geq 1$, $x\in\Z$. (Define $B_j(x)=\infty$ if $\bar{b}_x=0$.)
\item[--] $D_n(x)$ uniformly distributed on $[0,\bar{b}_x]$, $n\geq 1$, $x\in\Z$.
\item[--] $S_j(x)$ exponentially distributed with mean $1/\bar{c}_x$, $j\geq 1$, $x\in\Z$.
\item[--] $U_n^0(x)$ uniformly distributed on $[0,\bar{c}_x^0]$, $n\geq 1$, $x\in\Z$.
\item[--] $U_n^1(x)$ uniformly distributed on $[0,\bar{c}_x^1]$, $n\geq 1$, $x\in\Z$.
\end{itemize}
Moreover, for $n\geq 1$ and $x\in\Z$ define
\[
C_n(x)=\sum_{j=1}^n B_j(x) \quad \text{and} \quad
T_n(x)=\sum_{j=1}^n S_j(x).
\]
For a given initial configuration $\beta\in\{0,1\}^\Z$, define a process $\{\beta_t^\beta\}_{t\geq 0}$ from $\{C_n(x)\}$ and $\{D_n(x)\}$ as follows:
\begin{itemize}
\item[--] $\beta_0^\beta=\beta$,
\item[--] $\beta_s^\beta(x)$ flips from $0$ to $1$ iff $\beta_{s-}^\beta(x)=0$ and there exists an $n\geq 1$ such that $s=C_n(x)$ and $D_n(x)\geq \bar{b}_x-b(x,\beta_{s-}^\beta)$, 
\item[--] $\beta_s^\beta(x)$ flips from $1$ to $0$ iff $\beta_{s-}^\beta(x)=1$ and there exists an $n\geq 1$ such that $s=C_n(x)$ and $D_n(x)< b(x,\beta_{s-}^\beta)$.
\end{itemize}
By an approximation procedure, it is possible to prove that there exists a process with those properties and that such a process has flip rates $b(x,\beta)$. 

Given $\beta$,$\eta\in\{0,1\}^\Z$, we now define a process $\{\eta_t^{\beta,\eta}\}_{t\geq 0}$ from $\{\beta_t^\beta\}$, $\{T_n(x)\}$, $\{U_n^0(x)\}$ and $\{U_n^1(x)\}$ in the following way:
\begin{itemize}
\item[--] $\eta_0^{\beta,\eta}=\eta$,
\item[--] if $\beta_s^\beta(x)=0$, then $\eta_s^{\beta,\eta}(x)$ flips from $0$ to $1$ iff $\eta_{s-}^{\beta,\eta}(x)=0$ and there exists an $n\geq 1$ such that $s=T_n(x)$ and $U_n^0(x)\geq \bar{c}_x^0-\frac{\bar{c}_x^0}{ \bar{c}_x}c_0(x,\eta_{s-}^{\beta,\eta})$ and $\eta_s^{\beta,\eta}(x)$ flips from $1$ to $0$ iff $\eta_{s-}^{\beta,\eta}(x)=1$ and there exists an $n\geq 1$ such that $s=T_n(x)$ and $U_n^0(x)< \frac{\bar{c}_x^0}{\bar{c}_x}c_0(x,\eta_{s-}^{\beta,\eta})$,
\item[--] if $\beta_s^\beta(x)=1$, then $\eta_s^{\beta,\eta}(x)$ flips from $0$ to $1$ iff $\eta_{s-}^{\beta,\eta}(x)=0$ and there exists an $n\geq 1$ such that $s=T_n(x)$ and $U_n^1(x)\geq \bar{c}_x^1-\frac{\bar{c}_x^1}{\bar{c}_x}c_1(x,\eta_{s-}^{\beta,\eta})$ and $\eta_s^{\beta,\eta}(x)$ flips from $1$ to $0$ iff $\eta_{s-}^{\beta,\eta}(x)=1$ and there exists an $n\geq 1$ such that $s=T_n(x)$ and $U_n^1(x)< \frac{\bar{c}_x^1}{ \bar{c}_x}c_1(x,\eta_{s-}^{\beta,\eta})$.
\end{itemize}
It is clear that the process $\{(\beta_t^\beta,\eta_t^{\beta,\eta})\}$ has the correct flip rates. Moreover, the graphical representation gives us a coupling for all possible initial states and this coupling is exactly the maximal type coupling used in Section~\ref{paper2:proof_mainthm}. If we want to start the process at a random state with distribution $\rho$, we just add, independent of everything else, two random variables with joint distribution $\rho$. We then write $\{\beta_t^{\rho_1},\eta_t^{\rho_1,\rho_2}\}$ where $\rho_i$ denotes the $i$th marginal of $\rho$.

A possible proof of Conjecture~\ref{paper2:conj_pos_rates} may be based on the following lemma.
\begin{lemma}\label{paper2:ergodic_lemma}
If
\begin{equation}\label{paper2:ergodic_cond}
\liminf_{k\to\infty}\liminf_{t\to\infty}\P\left[\,\eta_t^{\beta,\emptyset}(x)=\eta_t^{\beta,\Z}(x),\,-k\leq x\leq k\,\right]>0
\end{equation}
for all $\beta\in\{0,1\}^\Z$, then $\nu_0=\nu_1$. 
\end{lemma}
\begin{proof}
From Lemma~\ref{paper2:statcoup} (or more precisely from the version of it with three processes) there exists a probability measure $\gamma$ on
\[
\left\{\,(\beta,\eta,\xi)\in U^3:\:\eta \leq  \xi\right\}
\]
which is stationary for $\{(\beta_t,\eta_t,\xi_t\}_{t\geq 0}$ and satisfies 
\[
\gamma_{12}=\nu_0, \quad \gamma_{13}=\nu_1 \quad\text{and}\quad\gamma_1=\mu,
\]
where $\mu$ is the unique stationary distribution for the background process. (Here, we use the same notation as in Lemma~\ref{paper2:statcoup}.) Our goal is to show that 
\[
\gamma\left(\eta=\xi\right)=1.
\] 
For given $k\geq 1$ and $t\geq 0$, we get
\begin{equation}\label{paper2:stateq1}
\begin{split}
&\gamma\left(\eta(x)=\xi(x),\,-k\leq x\leq k\right) 
=\gamma\left(\eta=\xi\right) \\
&\quad+\P\left[\,\eta_t^{\mu,\gamma_2}(x)=\eta_t^{\mu,\gamma_3}(x),\,-k\leq x\leq k\,|\,\eta_0^{\mu,\gamma_2}\not =\eta_0^{\mu,\gamma_3}\,\right]\left(1-\gamma\left(\eta=\xi\right)\right).
\end{split}
\end{equation}
Here, we have used that $\gamma$ is stationary and the fact that 
\[
\P\left[\,\eta_t^{\mu,\gamma_2}(x)=\eta_t^{\mu,\gamma_3}(x),\,-k\leq x\leq k\,|\,\eta_0^{\mu,\gamma_2}=\eta_0^{\mu,\gamma_3}\,\right]=1.
\]
From the inequalities
\[
\eta_t^{\mu,\emptyset}\leq\eta_t^{\mu,\gamma_2}\leq\eta_t^{\mu,\gamma_3}\leq\eta_t^{\mu,\Z},\quad t\geq 0,
\]
we get,
\begin{equation}\label{paper2:erg_ineq}
\begin{split}
&\P\left[\,\eta_t^{\mu,\gamma_2}(x)=\eta_t^{\mu,\gamma_3}(x),\,-k\leq x\leq k\,|\,\eta_0^{\mu,\gamma_2}\not =\eta_0^{\mu,\gamma_3}\,\right] \\
&\quad \geq \P\left[\,\eta_t^{\mu,\emptyset}(x)=\eta_t^{\mu,\Z}(x),\,-k\leq x\leq k\,|\,\eta_0^{\mu,\gamma_2}\not =\eta_0^{\mu,\gamma_3}\,\right].
\end{split}
\end{equation}
Moreover, from the graphical representation, we get that the events
\[
\{\,\eta_t^{\mu,\emptyset}(x)=\eta_t^{\mu,\Z}(x),\,-k\leq x\leq k\,\}\quad\text{and}\quad \{\,\eta_0^{\mu,\gamma_2}\not =\eta_0^{\mu,\gamma_3}\,\}
\]
are conditionally independent given the initial state of the background process and so we can write
\begin{equation}\label{paper2:cond_beta}
\begin{split}
&\P\left[\,\eta_t^{\mu,\emptyset}(x)=\eta_t^{\mu,\Z}(x),\,-k\leq x\leq k,\,\eta_0^{\mu,\gamma_2}\not =\eta_0^{\mu,\gamma_3}\,\right] \\
&\quad\int\P\left[\,\eta_t^{\mu,\emptyset}(x)=\eta_t^{\mu,\Z}(x),\,-k\leq x\leq k\,|\,\beta_0^\mu=\beta\,\right]\gamma\left(\eta\not =\xi\,|\,\beta\right)\,d\mu(\beta).
\end{split}
\end{equation}
Now, let us assume that
\[
\gamma\left(\eta\not =\xi\right)>0.
\]
Then
\[
\gamma\left(\eta\not =\xi\,|\,\beta\right)>0.
\]
on a set of positive $\mu$-measure. By using \eqref{paper2:ergodic_cond}, \eqref{paper2:cond_beta} together with Fatou's Lemma and then \eqref{paper2:erg_ineq}, we can conclude that
\[
\liminf_{k\to\infty}\liminf_{t\to\infty}\P\left[\,\eta_t^{\mu,\gamma_2}(x)=\eta_t^{\mu,\gamma_3}(x),\,-k\leq x\leq k\,|\,\eta_0^{\mu,\gamma_2}\not =\eta_0^{\mu,\gamma_3}\,\right]>0.
\]
However, by taking limits in \eqref{paper2:stateq1} we arrive at a contradiction and so we are done with the proof. 
\end{proof}
The question now is if it is possible to prove \eqref{paper2:ergodic_cond}. A natural first try is to fix the initial state of the background process and then proceed as in \cite[p.~393]{Gray_pos_rates} and define so called left and right edge processes. The properties on p.~394 and Proposition~2 on p.~395 are then easily verified. For the correlation property between the left and right edge processes, we can use \cite[Ch.~II, Corollary~2.21]{Liggett85} and since the Lemma in the proof of \cite[Theorem~1]{Gray_pos_rates} relies on the properties on \cite[p.~394]{Gray_pos_rates}, it may be possible to prove a version of it for our process. Having succeeded so far, there is some hard work left which we at the moment are not able to decide on if it is possible to do something similar or not. The only thing we can say is that the argument given in \cite[p.~399-403]{Gray_pos_rates} is based on a very similar construction as we have and if all the preliminary work go through, then there may be a quite good chance to get a full proof of Conjecture~\ref{paper2:conj_pos_rates}.   

\section*{Acknowledgment}

The author wants to thank Jeffrey Steif for valuable comments and especially for helping me with Lemma~\ref{paper2:statcoup}.
 

\bibliographystyle{amsplain}
\bibliography{/home/warfheimer/Backup_remote/Forskning/Paper1/references}



\end{document}